\documentclass[12pt]{amsart}

\newcommand{\bfG}{\mathbf G}
\newcommand{\bfB}{\mathbf B}
\newcommand{\bfH}{\mathbf H}
\newcommand{\bfT}{\mathbf T}
\newcommand{\bfU}{\mathbf U}

\newcommand{\bfc}{\mathbf c}

\DeclareMathOperator\GL{GL}
\DeclareMathOperator\U{U}

\newcommand\bZ{{\mathbb Z}}
\newcommand\bF{\mathbb F}

\newcommand\sub{\subseteq}

\usepackage{amssymb}
\usepackage{fullpage}
\usepackage{algorithmic}

\numberwithin{equation}{section}

\theoremstyle{plain}
\newtheorem{thm}[equation]{Theorem}

\theoremstyle{definition}
\newtheorem{exmp}[equation]{Example}

\theoremstyle{remark}
\newtheorem{rem}[equation]{Remark}

\title
{Conjugacy classes in Sylow $p$-subgroups of finite
Chevalley groups in bad characteristic}

\author{John D.~Bradley}
\address
{London School of Hygiene \& Tropical Medicine,
London, WC1E 7HT,
United Kingdom}
\email{john.bradley@lshtm.ac.uk}

\author{Simon M.~Goodwin}
\address{School of Mathematics,
University of Birmingham,
Birmingham, B15 2TT,
United Kingdom}
\email{s.m.goodwin@bham.ac.uk}


\date{\today}

\begin{document}

\begin{abstract}
Let $U = \bfU(q)$ be a Sylow $p$-subgroup of a finite Chevalley group $G = \bfG(q)$.
In \cite{GR} R\"ohrle and the second author determined a
parameterization of the conjugacy classes of $U$, for $\bfG$ of small
rank when $q$ is a power of a good prime for $\bfG$.  As a consequence they verified that the number
$k(U)$ of conjugacy classes of $U$ is given by a polynomial in $q$ with integer coefficients.
In the present paper, we consider the case when $p$ is a bad prime for $\bfG$.
We obtain a parameterization of the
conjugacy classes of $U$, when $\bfG$ has rank less than or equal to $4$, and $\bfG$ is not of type $F_4$.
In these cases we deduce that $k(U)$ is given by a polynomial in $q$ with integer coefficients;
this polynomial is different from the polynomial for good primes.
\end{abstract}

\maketitle

\section{Introduction}

Let $\bfG$ be a split simple algebraic group defined over the finite
field $\bF_p$, where $p$ is a prime.  For a closed subgroup $\bfH$ defined over
$\bF_p$ and a power $q$ of $p$, we write $H = \bfH(q)$ for the finite group
of $\bF_q$-rational points in $\bfH$.  Let $\bfB$ be a Borel subgroup
of $\bfG$ defined over $\bF_p$ and let $\bfU$ be the unipotent radical of $\bfG$.  Then
$U$ is a Sylow $p$-subgroup of the finite Chevalley group $G$.

In case $\bfG = \GL_n$, we may take $\bfU = \U_n$ to be the group of
upper unitriangular matrices.  In this situation, a longstanding conjecture
attributed to G.~Higman (cf.\ \cite{Hi}) says
that the number $k(\U_n(q))$ of conjugacy classes of $\U_n(q)$ is
as a function of $q$ a polynomial in $q$ with integer coefficients.
This has been verified for $n \le 13$ by A.~Vera-Lopez and
J.~M.~Arregi by computer calculation, see \cite{VA}; see also \cite{Ev}
for a recent alternative approach due to A.~Evseev.  Higman's conjecture
has attracted a great deal of interest, for example from G.~Robinson \cite{Ro}, J.~Thompson
\cite{Th} and J.~Alperin \cite{Al}.  Of course, Higman's conjecture
can be stated in terms of the irreducible complex characters of $\U_n(q)$,
where there has also been much interest, for example from G.~Lehrer \cite{Le} and
I.~M.~Isaacs \cite{Is}.

Recently, there has a lot of research in to the conjugacy classes and irreducible
characters of $U$ for arbitrary $\bfG$, see for example
\cite{Go1}, \cite{Go2}, \cite{GR}, \cite{HH}, \cite{HLM} and \cite{Sz}.
In particular,  G.~R\"ohrle and the second author computed (using GAP \cite{GAP})
a parameterization of the conjugacy classes of $U$, when $\bfG$ is simple of rank less than or equal than $6$ and $p$ is a good
prime for $\bfG$, see \cite{GR}.  In these
cases, the number $k(U)$ of conjugacy classes of $U$ is a
polynomial in $q$ with integer coefficients, which is independent of the prime $p$, \cite[Theorem 1.1]{GR}.

In this paper, we consider the conjugacy classes of $U$ when $p$
is a bad prime for $\bfG$.  We have developed an algorithm to determine a parameterization of the conjugacy classes of $U$ and implemented it
in the computer algebra system MAGMA \cite{MAGMA} using the functions
for groups of Lie type explained in \cite{CMT}.  Using the resulting computer program we have obtained a parameterization
when $\bfG$ has rank less than or equal to $4$, and $\bfG$ is not of type $F_4$.
As a consequence we have the following theorem.

\begin{thm}
Let $\bfG$ be a split simple algebraic group defined over the $\bF_p$ and let $\bfU$ be a maximal unipotent subgroup of G defined over $\bF_p$.
Suppose that $p$ is bad for $\bfG$, the rank of $G$ less than or equal to $4$, and $\bfG$ is not of type $F_4$.
Let $q$ be a power of $p$.
Then the number of conjugacy classes of $\bfU(q)$ is given by a
polynomial in $q$ with integer coefficients.
\end{thm}

Many of the results about the conjugacy classes in $U$ from \cite{Go1}
do not hold in bad characteristic, which means that the parameterization of the conjugacy classes
is more subtle.  As a particular consequence this means that the polynomial giving $k(U)$ for $p$ a bad prime is different from that for good primes.

The algorithm used for our calculations
is based on the algorithm from \cite{GR}, but it is adapted to calculate with
groups rather than Lie algebras, and allows us to deal with complications
that do not occur in good characteristic.  An explanation of the algorithm is given in
Section~\ref{S:expl}, then the results of our calculations are presented in Section~\ref{S:res}.
We remark that in turn the algorithm in \cite{GR} is based on the algorithm of H.~B\"urgstein
and W.~H.~Hesselink from \cite{BH}.  We recall that the aim in \cite{BH} was to
understand the adjoint $\bfB$-orbits in the Lie algebra of $\bfU$ for which
there is a great deal of motivation from geometric representation theory.

\section{Explanation of calculations} \label{S:expl}

\subsection{Notation}
Below we collect all the notation that we require to describe our algorithm.

Let $\bfG$ be a split simple algebraic group defined over $\bF_p$, where $p$ is a prime.  Let $\bfB$ be a Borel
subgroup of $\bfG$ defined over $\bF_p$ and let $\bfU$ be the unipotent radical of $\bfB$.  Fix a maximal torus $\bfT$ of $\bfG$
contained in $\bfB$ and defined over $\bF_p$, and let $\Phi$ be the root system of $\bfG$ with respect to $\bfT$.  Let $\Phi^+ \sub \Phi$ be the system
of positive roots determined by $\bfB$.  The partial
order on $\Phi$ determined by $\Phi^+$ is denoted by $\preceq$ and we let $N = |\Phi^+|$.

Let $q$ be a power of $p$.
Given a closed subgroup $\bfH$ of $\bfG$ that is defined over $\bF_p$,
we write $H = \bfH(q)$ for the group of $\bF_q$-rational points of $\bfH$.
Given $\alpha \in \Phi^+$ we let $U_\alpha$ be the corresponding
root subgroup of $U$, and we choose a parameterization $x_\alpha : \bF_q \to U_\alpha$.

We fix an enumeration of the set of positive roots $\Phi^+ = \{\beta_1,\dots,\beta_N\}$ such that $i \le j$ whenever $\beta_i
\preceq \beta_j$.   We abbreviate notation and write $x_i = x_{\alpha_i}$.
Define the sequence of normal subgroups
$$
U = M_0 \supseteq \dots \supseteq M_N = \{1\}
$$
of $U$ by
$$
M_i = \prod_{j=i+1}^N U_{\beta_j}.
$$
Let
$$
U_i = U/M_i
$$
and note that the conjugation action of $U$ on itself induces an
action of $U$ on $U_i$.  Given $yM_i \in U_i$, $x \in U$ and a subgroup $H$ of $U$, we write
$x \cdot yM_i = xyx^{-1}M_i$, we denote  the $H$-orbit of $yM_i$ by  $H \cdot yM_i$ and we write $C_U(yM_i)$ for the stabilizer of $yM_i$
in $U$.

\subsection{Algorithm}

We present an outline of the algorithm that we have used to
calculate a parameterization of the conjugacy
classes in $U$ when $p$ is a bad prime for $\bfG$ and the rank of $\bfG$ is small.
The idea of our algorithm
is to calculate the $U$-orbits in $U_i$ for successive $U_i$.  It is an adaptation of
the algorithm from \cite{GR}, which considers the case where $p$ is a good
prime and is based on results from \cite{Go1}.
The algorithm proceeds in steps to obtain a set of {\em minimal representatives} of the conjugacy classes in $U$.
We note that these representatives depend on our chosen enumeration of $\Phi^+$.
Before giving a formal outline of the algorithm we give a rough
description of how the algorithm works.

\subsubsection*{Rough description}

We work in steps.

\smallskip

\noindent
{\bf 0th step}:
We just note that $U_0 = \{M_0 = 1_{U_0}\}$, so that $M_0$ is the representative of the only $U$-orbit in $U_0$.
We set $Y_0 = \{M_0\}$ and $D(M_0) = U$.

\smallskip

\noindent
{\bf After $i$th step}: We have calculated a set $Y_i \sub U_i$ which is an approximation
to a parameterization of the $U$-orbits in $U_i$.  We have that $Y_i$
contains representatives of all the $U$-orbits in $U_i$ though some elements
of $Y_i$ may be representatives of the same $U$-orbit.  Also for each $yM_i \in Y_i$,
we have calculated an approximation $D(yM_i)$ of $C_U(yM_i)$.  We have that
$D(yM_i) \supseteq C_U(yM_i)$, but $D(yM_i)$ may be strictly larger than
$C_U(yM_i)$.  The key point is that if two elements of $yM_i, zM_i \in Y_i$ are $U$-conjugate,
then they are conjugate under $D(yM_i)$.
(We note that in the calculations that we have made, we
actually have $D(yM_i) = C_U(yM_i)$ for most $yM_i \in Y_i$, so that $Y_i$ is close
to a parameterization of the $U$-orbits in $U_i$.)

\smallskip

\noindent
{\bf $(i+1)$th step}: We consider each element $yM_i \in Y_i$ in turn.  We assume that
$y$ is of the form $y = \prod_{j=1}^i x_j(a_j)$, where $a_j \in \bF_q$.  Then we consider
$D(yM_i) \cdot yM_{i+1} \sub U_{i+1}$.

First suppose that $D(yM_i) = C_U(yM_i)$.  We have that
$$
C_U(yM_i) \cdot yM_{i+1} \sub \{y x_{i+1}(t) M_{i+1} \mid t \in \bF_q\}.
$$
We now consider three cases.

\noindent
{\em Ramification case}:
The easiest case to deal with is when we have $C_U(yM_i) \cdot yM_{i+1} = \{yM_{i+1}\}$.  Then we have that
no two elements of $\{y x_{i+1}(t) M_{i+1} \mid t \in \bF_q\}$ are conjugate under $C_U(yM_i)$.  Also we have
$C_U(y x_{i+1}(t)M_{i+1}) = C_U(yM_i)$ for all $t \in \bF_q$.  Thus we add $y x_{i+1}(t)M_{i+1}$ to $Y_{i+1}$ for each $t \in \bF_q$
and we set  $D(y x_{i+1}(t)M_{i+1}) = C_U(yM_i)$.

\noindent
{\em Inert case}:
The next case to consider is when we can see that $C_U(yM_i) \cdot yM_{i+1} = \{y x_{i+1}(t) M_{i+1} \mid t \in \bF_q\}$
by finding a morphism $x : \bF_q \to C_U(yM_i)$ such that $x(t) \cdot yM_{i+1} = y x_{i+1}(t) M_{i+1}$.  Using this morphism
we are able to calculate $C_U(yM_{i+1})$; we leave the details of how this is achieved until the formal
outline below.  Thus we add $yM_{i+1}$ to $Y_{i+1}$ along with $D(yM_{i+1}) = C_U(yM_{i+1})$.

\noindent
{\em Unresolved case}:
If neither the ramification case nor the inert case hold, then the algorithm does not try to calculate a parameterization of the
$U$-orbits with representatives in $\{y x_{i+1}(t) M_{i+1} \mid t \in \bF_q\}$.  We just add $y x_{i+1}(t)M_{i+1}$ to $Y_{i+1}$ for each $t \in \bF_q$ along with $D(y x_{i+1}(t)M_{i+1}) = C_U(yM_i)$.

\smallskip

When we have $D(yM_i) \supsetneqq C_U(yM_i)$ the situation is a little more complicated.  We have to consider
$(D(yM_i) \cdot yM_{i+1}) \cap \{y x_{i+1}(t) M_{i+1} \mid t_i \in \bF_q\}$ instead of $C_U(yM_i) \cdot yM_i$, and with
this modification we proceed to consider the three cases as above.

\smallskip

\noindent
{\bf After $N$th step}:  We end up with the set $Y_N$ and for each $y \in Y_N$ an approximation $D(y)$ of $C_U(y)$.
To complete our classification we have to calculate representatives of the $U$-orbits for which $D(y) \supsetneqq C_U(y)$.
This is achieved by considering these cases individually and making some calculations by hand.  We note that most of
these cases can be dealt with quite easily, and we give more details in the formal outline below, and also in
Example~\ref{ex:B4}.

\begin{rem}
We note that in case $p$ is good for $\bfG$, the results from \cite{Go1} tell us
that the unresolved case will not occur.  It is this possibility that makes the
situation for bad primes more complicated than that for good primes.
\end{rem}

\subsubsection*{An example}

Before attempting a formal outline of the algorithm we demonstrate it with an example.  We consider the
case where $\bfG$ is of type $C_3$ and $p=2$; we use the enumeration of $\Phi^+$ from MAGMA which is given in Table~\ref{Tab:C3roots} below.
We just consider those conjugacy classes of $U$ which have representatives of the form $y = x_2(a_2)x_3(a_3)y'$, where $a_2,a_3 \in \bF_q^\times$.
and $y' \in M_3$.  These correspond to the $2,3$ row of Table~\ref{Tab:C3} below.  This example is a lot less complicated than a number of cases that have
to be considered for $B_4$ and $C_4$, but it gives a good idea of how the algorithm works.

\begin{exmp} \label{ex:C3}
We easily see that after the 3rd step we have $yM_3 = x_2(a_2)x_3(a_3)M_3 \in Y_3$ for all $a_2,a_3 \in \bF_q^\times$, and that we have calculated $D(yM_3) = C_U(yM_3) = U$.
So elements of $C_U(yM_3)$ can be written in the form:
$$
x = x_1(t_1) \cdots x_9(t_9),
$$
where $t_1,\dots,t_9 \in \bF_q$.

To complete the 4th step we first calculate
$$
x \cdot yM_4 = x_2(a_2)x_3(a_3)x_4(a_2t_1)M_4.
$$
Therefore, we see that we are in the inert case, so only $yM_4 = x_2(a_2)x_3(a_3)M_4$ is added to $M_4$.  It is also easy to see that for
$x \in C_U(yM_4)$, we must have $t_1 = 0$. Thus we take $D(yM_4) = C_U(yM_4)$ and its elements are of the form:
$$
x = x_2(t_2) \cdots x_9(t_9),
$$
where $t_2,\dots,t_9 \in \bF_q$.

For the 5th step we consider
$$
x \cdot yM_5 = x_2(a_2)x_3(a_3)x_5(a_3t_2 + a_2t_3).
$$
So again we are in the inert case and we only add $yM_5 = x_2(a_2)x_3(a_3)M_5$ to $Y_5$.  To calculate $D(yM_5) = C_U(yM_5)$
we have the make the substitution $t_3 = \frac{a_3t_2}{a_2}$ in to the expression for $x$ above, so we see that
elements of $C_U(xM_5)$ are of the form:
$$
x = x_2(t_2) x_3\left(\frac{a_3t_2}{a_2}\right) x_4(t_4) \cdots x_9(t_9),
$$
where $t_2,t_4,\dots,t_9 \in \bF_q$.

For the 6th step we consider
$$
x \cdot yM_6 = x_2(a_2)x_3(a_3)x_6(a_3t_4)M_6.
$$
Again we are in the inert case, so just add $yM_6 = x_2(a_2)x_3(a_3)M_6$ to $Y_6$.  To determine $D(yM_6) = C_U(yM_6)$
we have to substitute $t_4 = 0$, and we get elements of the form
$$
x = x_2(t_2) x_3\left(\frac{a_3t_2}{a_2}\right) x_5(t_5) \cdots x_9(t_9),
$$
where $t_2,t_5,\dots,t_9 \in \bF_q$.

For the 7th step we consider
$$
x \cdot yM_7 = x_2(a_2)x_3(a_3)x_7(a_3t_2^2 + a_2 a_3t_2)M_7.
$$
We note that the map $t \mapsto a_3t^2 + a_2 a_3 t$ from $\bF_q$ to itself is not surjective.  Thus we are
in the unresolved case.  So we add $y = x_2(a_2) x_3(a_3) x_7(b_7) M_7$ to  $Y_7$ for all $b_7 \in \bF_q$, with
$D(x_2(a_2) x_3(a_3) x_7(b_7) M_7)$ consisting of elements of the form
$$
x = x_2(t_2) x_3\left(\frac{a_3t_2}{a_2}\right) x_5(t_5) \cdots x_9(t_9),
$$
where $t_2,t_5,\dots,t_9 \in \bF_q$.

For the 8th step we consider
$$
x \cdot yM_8 = x_2(a_2)x_3(a_3)x_7(a_3t_2^2 + a_2 a_3t_2 + b_7)x_8(a_2t_6)M_8.
$$
So we can see that we are in the inert case, so we just add $yM_8$ to $Y_8$.  Then to get
$D(yM_8)$, we have to substitute $t_6 = 0$, to get elements of the form
$$
x = x_2(t_2) x_3\left(\frac{a_3t_2}{a_2}\right) x_5(t_5) x_7(t_7) x_8(t_8) x_9(t_9).
$$

For the 9th and final step we consider
$$
x \cdot y = x_2(a_2)x_3(a_3)x_7(a_3t_2^2 + a_2 a_3t_2 + b_7).
$$
Thus we are in the ramification case, so we add $y =  x_2(a_2)x_3(a_3)x_7(b_7)x_9(b_9)$ to $Y_9$ for all
$b_9 \in \bF_q$.  We set $D(y) = D(yM_8)$.

The last thing we have to do is look at the one unresolved step.  We have
$$
x \cdot y = x_2(a_2)x_3(a_3)x_7(a_3t_2^2 + a_2 a_3t_2 + b_7)x_9(b_9).
$$
Now the map $\phi : t \mapsto a_3t^2 + a_2 a_3 t$ from $\bF_q$ to $\bF_q$ has kernel of size 2, so its image is
index 2 in $\bF_q$.  Therefore, we can use the action of $D(y_9)$ to assume that $b_7$ is either $0$ or
one chosen element that is not in the image of $\phi$.  We then write our elements as
$y = x_2(a_2)x_3(a_3)x_7(c_7)x_9(b_9)$, where $c_7$ is one of these two possibilities.  So we
see that in total we get $2(q-1)^2q$ conjugacy classes in this family, or $2v^2(v+1)$, where $v = q-1$.
\end{exmp}

\subsubsection*{Formal outline}

We now give our formal outline of how the algorithm works in general.
As suggested by the example, we actually group the representatives into families
depending on which $i$ there are nonzero values of $a_i$ in an expression of the
form $\prod_{i=1}^N x_i(a_i)$.  This is of course necessary in order to consider all different values of $q$ simultaneously.
Our outline of the algorithm is written in a sort of pseudo-code and we proceed in steps from the 0th step to the
$N$th step.

\smallskip
\noindent
{\bf Setup}:
We let $a_1,\dots,a_n$ and $t_1,\dots,t_N$ be indeterminates and consider the ring $R = \bF_p(a_1,\dots,a_N)[t_1,\dots,t_N]$.
We define $\mathcal G = \bfG(R)$ be the group of $R$-points of $G$, and we perform calculations in $\mathcal G$.

\smallskip
\noindent
{\bf 0th step}: We set:
\begin{itemize}
\item $C_0 = \{(\varnothing,\varnothing)\}$; and
\item $f_{(\varnothing,\varnothing),0,j}(a_k,t_l) = t_j \in R$ for $j = 1,\dots,N$.
\end{itemize}

\smallskip
\noindent
{\bf After $i$th step}:  We have defined:
\begin{itemize}
\item a subset $C_i$ of $\mathbb P(\{1,\dots,i\}) \times \mathbb P(\{1,\dots,i\})$ such that if $\bfc = (c,d) \in C_i$, then
$c \cap d = \varnothing$.
\item For each $\bfc \in C_i$ and $j = 1,\dots,N$, elements $f_{\bfc,i,j}(a_k,t_l) \in R$.
\end{itemize}
We explain how this corresponds to the rough description above.
Let $\bfc = (c,d) \in C_i$.  We set
$$
Y_\bfc = \{\prod_{j \in c \cup d} x_j(a_j)M_i \mid a_j \in \bF_q^\times \text{ for $j \in c$, } a_j \in \bF_q \text{ for $j \in d$, } \}.
$$
Then we have $Y_i = \bigcup_{\bfc \in C_i} Y_\bfc$.
Let $yM_i = \prod_{j \in c \cup d} x_j(a_j)M_i \in Y_\bfc$.
Then we have
$$
D(yM_i) = \{\prod_{j=1}^N x_j(f_{\bfc,j}(a_k,t_l)) \mid t_l \in \bF_q\}.
$$

\smallskip
\noindent
{\bf $(i+1)$th step}: \\
Initialize $C_{i+1} = \varnothing$.
\begin{algorithmic}
\FOR {$\bfc = (c,d) \in C_i$}
\STATE $y:= \prod_{j \in c \cup d} x_i(a_i)$
\STATE $x:= \prod_{j=1}^N x_i(f_{\bfc,j}(a_k,t_l))$
\STATE $zM_{i+1}:= x \cdot yM_{i+1}$
\STATE write $z = \prod_{j=1}^{i+1} x_j(g_j(a_k,t_l))$, where $g_j(a_k,t_l) \in R$
\STATE view $g_{i+1}(a_k,t_l)$ as a polynomial in $t_1,\dots,t_N$ with coefficients in $\bF_p(a_1,\dots,a_N)$
\IF[This is the ramification case]{$g_{i+1}(a_k,t_l) = 0$}
\STATE add $(c \cup \{i+1\},d)$ to $C_{i+1}$
\STATE set $f_{(c \cup \{i+1\},d),i+1,j}(a_k,t_l) = f_{\bfc,i,j}(a_k,t_l)$ for $j = 1,\dots,N$
\STATE add $\bfc$ to $C_{i+1}$
\STATE set $f_{\bfc,i+1,j}(a_k,t_l) := f_{\bfc,i,j}(a_k,t_l)$ for $j = 1,\dots,N$
\ELSE
\STATE Search for a variable $t_l$ such that $t_l$ occurs in $g_{i+1}(a_k,t_l)$ only in a single term of the form
$h(a_k)t_l$, where $h(a_k)$ is a monomial in the $a_k$, (i.e.\ it is of the form $h(a_k) = \frac{a_1^{m_1} \dots a_N^{m_N}}{a_1^{n_1} \dots a_N^{n_N}}$, where
$m_1,\dots,n_N,n_1,\dots,n_N \in \bZ_{\ge 0}$), and $t_l$ does not occur in $g_j(a_k,t_l)$ for any $j \in d$.
  \IF[This is the inert case]{such a variable $t_l$ has been found}
  \STATE choose $l$ to be maximal
  \STATE add $\bfc$ to $C_{i+1}$
  \STATE define $f_{\bfc,i+1,j}(a_k,t_l)$ by substituting $t_l:=  \frac{g_{i+1}(a_k,t_l) - h(a_k)t_l}{h(a_k)}$ in to $f_{\bfc,i,j}(a_k,t_l)$ for $j = 1,\dots,N$
  \ELSE[This is the unresolved case]
  \STATE add $(c,d \cup \{i+1\})$ to $C_{i+1}$
  \STATE set $f_{(c,d \cup \{i+1\}),i+1,j}(a_k,t_l) := f_{\bfc,i,j}(a_k,t_l)$
  \ENDIF
\ENDIF
\ENDFOR
\end{algorithmic}

\smallskip
\noindent
{\bf After $N$th step}:
Output the set $C_N$ and for each $\bfc = (c,d) \in C_N$, the polynomials $f_{\bfc,N,j}(a_k,t_l)$ for which $j \in d$.

We then have to analyze this output to determine the parameterization of the conjugacy classes of $U$ and to calculate
$k(U)$.  To do this for $\bfc = (c,d) \in C_N$, $j \in d$ and fixed $a_k \in \bF_q$, we consider  $f_{\bfc,N,j}(a_k,t_l)$ as a function from $\bF_q^N$ to $\bF_q^N$ in the
$t_l$.  We consider the image of this function and this allows us to make an assumption on the value of $a_j$ in a representative of
the conjugacy class of $U$ of the form $y = \prod_{i \in c \cup d} x_i(a_i)$.  This method is best illustrated by an example given below.

\smallskip

Before giving an example of how the output is analyzed, we mention a simplification to our algorithm
obtained by using the conjugation action of $T$.
Let $\bfc = (c,d) \in C_i$ for some $i$.  Let $c' \sub c$ and suppose that for any $y= \prod_{j \in c \cup d} x_j(a_j)$
there is $t \in T$ such that $tyt^{-1} = y' = \prod_{j \in c \cup d} x_j(a_j')$, where $a'_j = 1$ for $j \in c'$.
Then, since $C_U(y) = tC_U(y')t^{-1}$, the rational functions
$f_{\bfc,i,j}(a_k,t_l), f_{\bfc,i,j}(a'_k,t_l) \in R$ are related ``via conjugation by $t$''.
This implies that in determining which of the inert, ramification or unresolved case holds in
later steps, we can just look at elements of the form $y= \prod_{j \in c \cup d} x_j(a_j)$, where $a_j = 1$ for $j \in c'$, i.e.\ we
can normalize these coefficients to 1.  This simplification reduces the number of indeterminates $a_j$ that we need to involve,
which speeds up the algorithm considerably.  We just need to note that for each representative of a conjugacy class of $U$ of the
form $y= \prod_{j \in c \cup d} x_j(a_j)$, where $a_j = 1$ for $j \in c'$ that
we obtain accounts for $(q-1)^{|c'|}$ conjugacy classes of $U$.

To check whether the above condition holds for $c'$, we just need to calculate
the Smith normal form of the integer matrix with rows given by the expressions
for $\{\beta_j \mid j \in c'\}$ in terms of the simple roots.  The condition
holds if this Smith normal form has one entry equal to 1 on each row.

\smallskip

In our example below of how we analyze the output, we consider output obtained from the algorithm when
$\bfG$ is of type $B_4$, and
the algorithm had normalized certain coefficients to $1$ using the action of $T$.

\begin{exmp} \label{ex:B4}
We look at a $B_4$ case, $p=2$, where the output gives $\bfc = (c,d)$, where $|c| = 4$  and $|d| = 3$.  We have the following
$f_{\bfc,N,j}(a_k,t_l)$ for $j \in d$.
\begin{itemize}
\item $(a_{10}+1)t_5$;
\item $t_2^2 + t_2 + t_6^2 +a_{10}t_6$; and
\item $(a_{10}+1)t_5$.
\end{itemize}
To analyze this, we have to consider the cases $a_{10} = 1$ and $a_{10} \ne 1$ separately.

First suppose that $a_{10} \ne 1$.  Then we see that the first function is surjective, so we can assume that
the corresponding $a_i=0$.  Next we look at
$$
t_2^2 + t_2 + t_6^2 +a_{10}t_6 = (t_2+t_6)^2 + (t_2+t_6) + (a_{10}+1)t_6,
$$
which is seen to be surjective by taking $t_2 = t_6$ so we get $(a_{10}+1)t_6$.  Thus again
we can assume that the corresponding $a_i = 0$.
Finally, we look at the last function which is equal to first, but this can now longer be used to
make an assumption on the corresponding $a_i$ as it must be zero to centralize the earlier one.  So we see that
we get $v^3(v-1)(v+1)$ classes, where $v = q-1$.

In case $a_{10} = 1$, the first and last function are zero so we can make no assumptions on the corresponding $a_i$'s.  The middle
function is $(t_2+t_6)^2 + (t_2+t_6)$, which has the same image as the map $t \mapsto t^2 + t$ from $\bF_q$ to itself, which is
of size $q/2$.  So as in the 7th step in Example~\ref{ex:C3} we can assume that the corresponding $a_i$ is one of two elements of
$\bF_q$.  So we get $2v^3(v+1)^2$ classes.

Therefore, in total we get $3v^5+4v^4+v^3$ conjugacy classes of $U$ in this family.
\end{exmp}

\subsection{Implementation in MAGMA}

We now briefly describe how we have programmed our algorithm in MAGMA \cite{MAGMA}.  This is done using the
functions in MAGMA for groups of Lie Type explained in \cite{CMT}.  We can use the following code to define the group $\mathcal G = \bfG(R)$
used for calculations in our algorithm:
\begin{verbatim}
F:=RationalFunctionField(FiniteField(p),24);
R:=PolynomialRing(F,24);
G:=GroupOfLieType(type,R : Normalising:=true, Method:="CollectionFromLeft");
\end{verbatim}
The extra arguments in the definition of {\tt G} are required to ensure that MAGMA orders elements
$y:=\prod_{i=1}^N x_i(a_i)$ in ascending order.  Once we have defined {\tt G} we can carry out
all the operations required in our algorithm in MAGMA.

\section{Results} \label{S:res}

In this section we present the results of our computations.  First in Table~\ref{Tab:polys} below, we give the polynomials
giving $k(U)$.  We include both the polynomials for good and bad characteristics, so that the difference can be seen.  We write these
polynomials in $v = q-1$, and observe that in the case of bad characteristic the coefficients are positive,
as is the case for good characteristic.   As is the
case for good characteristic with have observed that the polynomials giving $k(U)$ for $\bfG$ of type $B_l$ and $C_l$ are
the same for $l \le 4$, so we only record these polynomials once.  It would be very interesting to have an explanation for this
phenomenon -- the representatives for the conjugacy classes of $U$ given in Tables~\ref{Tab:B3} and \ref{Tab:C3} suggest
that there is not a natural bijection.

\begin{table}[h!tb]
\renewcommand{\arraystretch}{1.4}
\begin{tabular}{|l|l|l|}
\hline
Type & Prime & $k(U(q))$  \\
\hline\hline
$B_2$ & $\ne 2$ & $2v^2 + 4v + 1$ \\
\cline{2-3}
& $= 2$ & $5v^2+4v+1$  \\
\hline
$G_2$ & $\ne 2,3$ & $2v^2 + 4v + 1$ \\
\cline{2-3}
& $= 3$ & $2v^3+11v^2+6v+1$ \\
\cline{2-3}
& $= 2$ & $v^3+8v^2+6v+1$ \\
\hline
$B_3$, $C_3$ & $\ne 2$ & $v^4+8v^3+16v^2+9v+1$ \\
\cline{2-3}
& $= 2$ & $2v^4+19v^3+25v^2+9v+1 $  \\
\hline
$B_4, C_4$ & $\ne 2$ & $v^6+11v^5+48v^4+88v^3+64v^2+16v+1$ \\
\cline{2-3}
& $= 2$ & $2v^6+31v^5+136v^4+168v^3+82v^2+16v+1$  \\
\hline
$D_4$ & $\ne 2$ & $2v^5+15v^4+36v^3+34v^2+12v+1$ \\
\cline{2-3}
& $= 2$ & $2v^5+18v^4+36v^3+34v^2+12v+1$  \\
\hline
\end{tabular}
\medskip
\caption{Polynomials in $v = q-1$ giving $k(U)$} \label{Tab:polys}
\end{table}

In the even numbered tables below we give the parameterization of the conjugacy classes of  $U$ when the rank of
$\bfG$ is 3 or less.  We do not include such tables for the rank 4 cases, as these would take up a lot of space; it is
straightforward to generate this information from the computer programme.  In the odd numbered tables, we give the enumeration of the positive roots
used for the parameterization of the conjugacy classes of $U$.  These roots are given
as expressions in terms of the simple roots as labelled in \cite[Planches I-IX]{Bo}.

For the even numbered tables, the conjugacy classes are put in to families, and given a name depending on the nonzero
entries in the representatives.  The families are given both for good and bad primes, so that the difference
can be seen.  We use the following notation to describe the elements of the families.
\begin{itemize}
\item $a_i$ is any element of $\bF_q^\times$;
\item $b_i$ is any element of $\bF_q$;
\item $c_i$ is one of two possibilities if $p=2$ and one of three possibilities if $p=3$.
\begin{itemize}
\item For $p=2$ it is either $0$ or an element of $\bF_q$ in the nonzero coset of the image of a map
from $\bF_q$ to itself of the form $t \to  \alpha t^2 +\beta t$, where $\alpha,\beta \in \bF_q^\times$, as in Example~\ref{ex:C3}.
\item For $p=3$ it is either $0$ or one of the two elements of $\bF_q$ in the nonzero cosets of the image of a map
from $\bF_q$ to itself of the form $t \to \alpha t^3 + \beta t$, where $\alpha,\beta \in \bF_q^\times$.
\end{itemize}
In both cases the elements $\alpha, \beta \in \bF_q^\times$ depend on the $a_j$ for $j < i$, and, in fact,
are rational functions in the $a_j$.
\item $d_i$ in any element of $\bF_q^\times \setminus \{\alpha\}$, where $\alpha$ depends on the $a_j$ for $j < i$ as a rational function.
\item $e_i$ is a single element of $\bF_q^\times$, which depends on the $a_j$ for $j < i$ as a rational function.
\item $f_i$ is either 0 or a single element of $\bF_q^\times$, which depends on the $a_j$ for $j < i$ as a rational function.
\end{itemize}
In the fourth column we give the size of the family as a polynomial in $v = q-1$.  The families are chosen
so that the size of the centralizer of elements of a family are equal, in the last column we give this centralizer size.

\begin{table}[h!tb]
\renewcommand{\arraystretch}{1.4}
\begin{tabular}{|l|l|l|l|l|}
\hline
Name & Prime & Family & Size of family & Centralizer size  \\
\hline\hline
1,2 & $\ne 2$ & $x_1(a_1)x_2(a_2)$ & $v^2$ & $q^2$ \\
\cline{2-5}
& $= 2$ & $x_1(a_1)x_2(a_2)x_4(c_4)$ & $2v^2$ & $2q^2$ \\
\hline
1 &- & $x_1(a_1)x_4(b_4)$ & $v(v+1)$ & $q^3$ \\
\hline
2 & $\ne 2$ & $x_2(a_2)$  & $v$ & $q^2$ \\
\cline{2-5}
& $= 2$ & $x_2(a_2)x_4(b_4)$  & $v(v+1)$ & $q^3$ \\
\hline
3 & $\ne 2$ & $x_3(a_3)$  & $v$ & $q^3$ \\
\cline{2-5}
& $= 2$ & $x_3(a_3)x_4(b_4)$  & $v(v+1)$ & $q^4$ \\
\hline
4 & - & $x_4(b_4)$  & $v+1$ & $q^4$ \\
\hline
\end{tabular}
\medskip
\caption{Conjugacy classes of $U$ for $G$ of type $B_2$} \label{Tab:B2}
\end{table}

\begin{table}[h!tb]
\renewcommand{\arraystretch}{1.4}
\begin{tabular}{|l|l|l|l|l|l|l|l|}
\hline
1 & \begin{tabular}{ll} 1 & 0 \\ \end{tabular}
  &
2 & \begin{tabular}{ll} 0 & 1 \\ \end{tabular}
  &
3 & \begin{tabular}{ll} 1 & 1 \\ \end{tabular}
  &
4 & \begin{tabular}{ll} 1 & 2 \\ \end{tabular}
  \\ \hline
\end{tabular}
\medskip
\caption{Enumeration of positive roots for $B_2$} \label{Tab:B2roots}
\end{table}

\newpage

\begin{table}[h!tb]
\renewcommand{\arraystretch}{1.4}
\begin{tabular}{|l|l|l|l|l|}
\hline
Name & Prime & Family & Size of family & Centralizer size  \\
\hline\hline
 1,2 & $\ne 2,3$ & $x_1(a_1)x_2(a_2)$ & $v^2$ & $q^2$ \\
\cline{2-5}
& $= 3$ & $x_1(a_1)x_2(a_2)x_5(c_5)$ & $3v^2$ & $3q^2$ \\
\cline{2-5}
& $= 2$ & $x_1(a_1)x_2(a_2)x_4(c_4)$ & $2v^2$ & $2q^2$ \\
\hline
1 & $\ne 2,3$ & $x_1(a_1)x_6(b_6)$ & $v(v+1)$ & $q^3$ \\
\cline{2-5}
& $= 3$ & $x_1(a_1)x_5(b_5)x_6(b_6)$ & $v(v+1)^2$ & $q^4$ \\
\cline{2-5}
& $= 2$ & $x_1(a_1)x_4(a_4)x_6(c_6)$ & $2v^2$ & $2q^4$ \\
&  & $x_1(a_1)$ & $v$ & $q^3$ \\
\hline
2 & $-$ & $x_2(a_2)x_4(b_4)x_5(b_5)$ & $v(v+1)^2$ & $q^4$ \\
\hline
3 & $\ne 2,3$ & $x_3(a_3)x_5(b_5)$ & $v(v+1)$ & $q^4$ \\
\cline{2-5}
& $= 3$ & $x_3(a_3)x_5(a_5)$ & $v^2$ & $q^4$ \\
&  & $x_3(a_3)x_6(b_6)$ & $v(v+1)$ & $q^5$ \\
\cline{2-5}
& $= 2$ & $x_3(a_3)x_4(a_4)x_5(c_5)$ & $2v^2$ & $2q^4$ \\
&  & $x_3(a_3)$ & $v$ & $q^4$ \\
\hline
4, & $\ne 3$ & $x_4(a_4)$ & $v$ & $q^4$ \\
\cline{2-5}
&$= 3$ & $x_4(a_4)x_5(a_5)$ & $v^2$ & $q^5$ \\
& & $x_4(a_4)x_5(b_6)$ & $v(v+1)$ & $q^6$ \\
\hline
5 & $-$ & $x_5(a_5)$ & $v$ & $q^5$ \\
\hline
6 & $\ne 2,3$ & $x_6(b_6)$ & $v+1$ & $q^6$ \\
\hline
\end{tabular}
\medskip
\caption{Conjugacy classes of $U$ for $G$ of type $G_2$}
\label{Tab:G2}
\end{table}

\begin{table}[h!tb]
\renewcommand{\arraystretch}{1.4}
\begin{tabular}{|l|l|l|l|l|l|}
\hline
1 & \begin{tabular}{ll} 1 & 0\end{tabular} &
2 & \begin{tabular}{ll} 0 & 1\end{tabular} &
3 & \begin{tabular}{ll} 1 & 1\end{tabular} \\ \hline
4 & \begin{tabular}{ll} 2 & 1\end{tabular} &
5 & \begin{tabular}{ll} 3 & 1\end{tabular} &
6 & \begin{tabular}{ll} 3 & 2\end{tabular} \\ \hline
\end{tabular}
\medskip
\caption{Enumeration of positive roots for $G_2$} \label{Tab:G2roots}
\end{table}

\bigskip

{\em Acknowledgments}: We would like to thank the Universidad de los Andes
for supporting a research visit of the second author during which some of this research was completed.  We also thank Scott
Murray for providing helpful answers to some questions about MAGMA.

\newpage

\tiny

\begin{samepage}

\begin{table}[h!tb]
\renewcommand{\arraystretch}{1.1}
\begin{tabular}{|l|l|l|l|l|}
\hline
Name & Prime & Family & Size of family & Centralizer size  \\
\hline\hline
1,2,3 & $\ne 2$ & $x_1(a_1)x_2(a_2)x_3(a_3)$ & $v^3$ & $q^3$ \\
\cline{2-5}
      & $= 2$ & $x_1(a_1)x_2(a_2)x_3(a_3)x_7(c_7)$ & $2v^3$ & $2q^3$ \\
\hline
1,2 & - & $x_1(a_1)x_2(a_2)x_7(b_7)$ & $v^2(v+1)$ & $q^4$ \\
\hline
1,3,5 & $\ne 2$ & $x_1(a_1)x_3(a_3)x_5(a_5)$ & $v^3$ & $q^4$ \\
\cline{2-5}
      & $= 2$ & $x_1(a_1)x_3(a_3)x_5(a_5)x_7(d_7)$ & $v^3(v-1)$ & $q^5$ \\
      &       & $x_1(a_1)x_3(a_3)x_5(a_5)x_7(e_7)x_9(c_9)$ & $2v^3$ & $2q^5$ \\
      &       & $x_1(a_1)x_3(a_3)x_5(a_5)$ & $v^3$ & $q^5$ \\
\hline
1,3 & $\ne 2$ & $x_1(a_1)x_3(a_3)x_9(b_9)$ & $v^2(v+1)$ & $q^5$ \\
\cline{2-5}
    &   $= 2$ & $x_1(a_1)x_3(a_3)x_7(a_7)x_9(c_9)$ & $2v^3$ & $2q^6$ \\
    &         & $x_1(a_1)x_3(a_3)$ & $v^2$ & $q^6$ \\
\hline
1,5 & $\ne 2$ & $x_1(a_1)x_5(a_5)$ & $v^2$ & $q^4$ \\
\cline{2-5}
    & $= 2$ & $x_1(a_1)x_5(a_5)x_7(a_7)$ & $v^3$ & $q^5$ \\
    &     & $x_1(a_1)x_5(a_5)x_9(c_9)$ & $2v^2$ & $2q^5$ \\
\hline
1,7 & - & $x_1(a_1)x_5(a_7)$ & $v^2$ & $q^5$ \\
\hline
1 & - & $x_1(a_1)x_9(b_9)$ & $v(v+1)$ & $q^6$ \\
\hline
2,3 & $\ne 2$ & $x_2(a_2)x_3(a_3)$ & $v^2$ & $q^3$ \\
\cline{2-5}
  & $= 2$ & $x_2(a_2)x_3(a_3)x_7(c_7)x_8(b_8)$ & $2v^2(v+1)$ & $2q^4$ \\
\hline
2 & - & $x_2(a_2)x_6(b_6)x_7(b_7)x_8(b_8)$ & $v(v+1)^3$ & $q^6$ \\
\hline
3,4 & $\ne 2$ & $x_3(a_3)x_4(a_4)$ & $v^2$ & $q^4$ \\
\cline{2-5}
    & $= 2$ & $x_3(a_3)x_4(a_4)x_7(a_7)$ & $v^3$ & $q^5$ \\
    &       & $x_3(a_3)x_4(a_4)x_8(c_8)$ & $2v^2$ & $2q^5$ \\
\hline
3 & $\ne 2$ & $x_3(a_3)x_9(b_9)$ & $v(v+1)$ & $q^5$ \\
\cline{2-5}
   & $=2$      & $x_3(a_3)x_7(a_7)x_9(b_9)$ & $v^2(v+1)$ & $q^6$ \\
  &       & $x_3(a_3)x_8(b_8)x_9(b_9)$ & $v(v+1)^2$ & $q^7$ \\
\hline
4,5 & $\ne 2$ & $x_4(a_4)x_5(a_5)x_8(b_8)$ & $v^2(v+1)$ & $q^6$ \\
\cline{2-5}
  & $= 2$ & $x_4(a_4)x_5(a_5)x_7(a_7)x_8(c_8)$ & $2v^3$ & $2q^6$ \\
  &       & $x_4(a_4)x_5(a_5)$ & $v^2$ & $q^6$ \\
\hline
4,7 & - & $x_4(a_4)x_7(a_7)$ & $v^2$ & $q^6$ \\
\hline
4 & - & $x_4(a_4)x_8(b_8)$ & $v(v+1)$ & $q^7$ \\
\hline
5 & $\ne 2$ & $x_5(a_5)x_8(b_8)$ & $v(v+1)$ & $q^6$ \\
\cline{2-5}
& $= 2$ & $x_5(a_5)x_7(a_7)x_8(b_8)$ & $v^2(v+1)$ & $q^7$ \\
&  & $x_5(a_5)x_8(a_8)$ & $v^2$ & $q^7$ \\
&  & $x_5(a_5)x_9(b_9)$ & $v(v+1)$ & $q^8$ \\
\hline
6 & $\ne 2$ & $x_6(a_6)x_7(b_7)$ & $v(v+1)$ & $q^7$ \\
\cline{2-5}
& $= 2$ & $x_6(a_6)x_7(a_7)$ & $v^2$ & $q^7$ \\
&  & $x_6(a_6)x_8(a_8)$ & $v^2$ & $q^8$ \\
&  & $x_6(a_6)x_9(b_9)$ & $v(v+1)$ & $q^9$ \\
\hline
7 & - & $x_7(a_7)$ & $v$ & $q^7$ \\
\hline
8 & - & $x_8(a_8)$ & $v$ & $q^8$ \\
\hline
9 & - & $x_9(b_9)$ & $v+1$ & $q^9$ \\
\hline
\end{tabular}
\medskip
\caption{Conjugacy classes of $U$ for $G$ of type $B_3$}
\label{Tab:B3}
\end{table}

\begin{table}[h!tb]
\renewcommand{\arraystretch}{1.1}
\begin{tabular}{|l|l|l|l|l|l|}
\hline
1 & \begin{tabular}{lll} 1 & 0 & 0 \\ \end{tabular}
  &
2 & \begin{tabular}{lll} 0 & 1 & 0 \\ \end{tabular}
  &
3 & \begin{tabular}{lll} 0 & 0 & 1 \\ \end{tabular}
  \\ \hline
4 & \begin{tabular}{lll} 1 & 1 & 0 \\ \end{tabular}
  &
5 & \begin{tabular}{lll} 0 & 1 & 1 \\ \end{tabular}
  &
6 & \begin{tabular}{lll} 1 & 1 & 1 \\ \end{tabular}
  \\ \hline
7 & \begin{tabular}{lll} 0 & 1 & 2 \\ \end{tabular}
  &
8 & \begin{tabular}{lll} 1 & 1 & 2 \\ \end{tabular}
  &
9 & \begin{tabular}{lll} 1 & 2 & 2 \\ \end{tabular}
  \\ \hline

\end{tabular}
\medskip
\caption{Enumeration of positive roots for $B_3$} \label{Tab:B3roots}
\end{table}

\end{samepage}

\newpage

\begin{samepage}

\begin{table}[h!tb]
\renewcommand{\arraystretch}{1.1}
\begin{tabular}{|l|l|l|l|l|}
\hline
Name & Prime & Family & Size of family & Centralizer size  \\
\hline\hline
1,2,3 & $\ne 2$ & $x_1(a_1)x_2(a_2)x_3(a_3)$ & $v^3$ & $q^3$ \\
\cline{2-5}
& $=2$ & $x_1(a_1)x_2(a_2)x_3(a_3)x_7(c_7)$ & $2v^3$ & $2q^3$ \\
\hline
1,2 &- & $x_1(a_1)x_2(a_2)x_7(b_7)$ & $v^2(v+1)$ & $q^4$ \\
\hline
1 & $\ne 2$ & $x_1(a_1)x_3(b_3)x_5(b_5)x_7(b_7)$ & $v(v+1)^3$ & $q^5$ \\
\cline{2-5}
& $= 2$ & $x_1(a_1)x_3(a_3)x_5(a_5)x_7(d_7)$ & $v^3(v-1)$ & $q^5$ \\
  &  & $x_1(a_1)x_3(a_3)x_5(a_5)x_7(f_7)x_9(c_9)$ & $4v^3$ & $2q^5$ \\
  & & $x_1(a_1)x_3(a_3)x_7(b_7)$ & $v^2(v+1)$ & $q^5$ \\
  & & $x_1(a_1)x_5(a_5)x_7(b_7)$ & $v^2(v+1)$ & $q^5$ \\
  & & $x_1(a_1)x_7(a_7)x_9(c_9)$ & $2v^2$ & $2q^5$ \\
  & & $x_1(a_1)x_9(b_9)$ & $v(v+1)$ & $q^6$ \\
\hline
2,3 & $\ne 2$ & $x_2(a_2)x_3(a_3)x_9(b_9)$ & $v^2(v+1)$ & $q^4$ \\
\cline{2-5}
    & $= 2$ & $x_2(a_2)x_3(a_3)x_7(c_7)x_9(b_9)$ & $2v^2(v+1)$ & $q^4$ \\
\hline
2,6 & $\ne 2$ & $x_2(a_2)x_6(a_6)$ & $v^2$ & $q^4$ \\
\cline{2-5}
& $= 2$ & $x_2(a_2)x_6(a_6)x_7(b_7)x_9(b_9)$ & $v^2(v+1)^2$ & $q^6$ \\
\hline
2 & $\ne 2$ & $x_2(a_2)x_9(b_9)$ & $v(v+1)$ & $q^5$ \\
\cline{2-5}
  & $= 2$ & $x_2(a_2)x_7(b_7)x_9(b_9)$ & $v(v+1)^2$ & $q^6$ \\
\hline
3,4 & $\ne 2$ & $x_3(a_3)x_4(a_4)x_7(b_7)$ & $v^2(v+1)$ & $q^5$ \\
\cline{2-5}
  & $= 2$ & $x_3(a_3)x_4(a_4)x_7(a_7)$ & $v^3$ & $q^5$ \\
    &  & $x_3(a_3)x_4(a_4)x_9(c_9)$ & $2v^2$ & $q^5$ \\
\hline
3,7 & - & $x_3(a_3)x_7(a_7)x_9(b_9)$ & $v^2(v+1)$ & $q^5$ \\
\hline
3,8 & $\ne 2$ & $x_3(a_3)x_8(a_8)$ & $v^2$ & $q^6$ \\
\cline{2-5}
    & $= 2$ & $x_3(a_3)x_8(a_8)x_9(b_9)$ & $v^2(v+1)$ & $q^5$ \\
\hline
3 &  & $x_3(a_3)x_9(b_9)$ & $v(v+1)$ & $q^7$ \\
\hline
4,5 & $\ne 2$ & $x_4(a_4)x_5(a_5)$ & $v^2$ & $q^5$ \\
\cline{2-5}
    & $= 2$ & $x_4(a_4)x_5(a_5)x_7(a_7)x_9(c_9)$ & $2v^3$ & $2q^6$ \\
    & & $x_4(a_4)x_5(a_5)$ & $v^2$ & $q^6$ \\
\hline
4 & $\ne 2$ & $x_4(a_4)x_7(b_7)$ & $v(v+1)$ & $q^6$ \\
\cline{2-5}
  & $= 2$ & $x_4(a_4)x_7(a_7)$ & $v^2$ & $q^6$ \\
  &       & $x_4(a_4)x_9(b_9)$ & $v(v+1)$ & $q^7$ \\
\hline
5 & $\ne 2$ & $x_5(a_5)x_9(b_9)$ & $v(v+1)$ & $q^6$ \\
\cline{2-5}
  & $= 2$ & $x_5(a_5)x_7(b_7)x_9(b_9)$ & $v(v+1)^2$ & $q^7$ \\
\hline
6 & $\ne 2$ & $x_6(a_6)x_7(b_7)$ & $v(v+1)$ & $q^7$ \\
\cline{2-5}  & $= 2$ & $x_6(a_6)x_7(a_7)$ & $v^2$ & $q^7$ \\
  & & $x_6(a_6)x_9(b_9)$ & $v(v+1)$ & $q^7$ \\
\hline
7 & - & $x_7(a_7)x_9(b_9)$ & $v(v+1)$ & $q^8$ \\
\hline
8 & $\ne 2$ & $x_8(a_8)$ & $v$ & $q^8$ \\
\cline{2-5}  & $= 2$ & $x_8(a_8)x_9(b_9)$ & $v(v+1)$ & $q^9$ \\
\hline
9 & - & $x_9(b_9)$ & $v+1$ & $q^9$ \\
\hline
\end{tabular}
\medskip
\caption{Conjugacy classes of $U$ for $G$ of type $C_3$}
\label{Tab:C3}
\end{table}

\begin{table}[h!tb]
\renewcommand{\arraystretch}{1.1}
\begin{tabular}{|l|l|l|l|l|l|}
\hline
1 & \begin{tabular}{lll} 1 & 0 & 0 \\ \end{tabular}
  &
2 & \begin{tabular}{lll} 0 & 1 & 0 \\ \end{tabular}
  &
3 & \begin{tabular}{lll} 0 & 0 & 1 \\ \end{tabular}
  \\ \hline
4 & \begin{tabular}{lll} 1 & 1 & 0 \\ \end{tabular}
  &
5 & \begin{tabular}{lll} 0 & 1 & 1 \\ \end{tabular}
  &
6 & \begin{tabular}{lll} 1 & 1 & 1 \\ \end{tabular}
  \\ \hline
7 & \begin{tabular}{lll} 0 & 2 & 1 \\ \end{tabular}
  &
8 & \begin{tabular}{lll} 1 & 2 & 1 \\ \end{tabular}
  &
9 & \begin{tabular}{lll} 2 & 2 & 1 \\ \end{tabular}
  \\ \hline
\end{tabular}
\medskip
\caption{Enumeration of positive roots for $C_3$}
\label{Tab:C3roots}
\end{table}

\end{samepage}

\end{document}